\documentclass[10pt]{article}
\usepackage{amsmath,amssymb,amsthm}
\usepackage{hyperref}
\usepackage{units}
\usepackage{color}
\usepackage[T1]{fontenc}
\usepackage[utf8]{inputenc}
\usepackage{authblk}
\usepackage{bm,latexsym,mathrsfs,enumerate}
\DeclareMathOperator{\Li}{Li}

\title{New Finite and Infinite Summation Identities Involving the Generalized Harmonic Numbers\thanks{%
MSC 2010: 65B10, 11B99}}
\author[1]{Kunle Adegoke\thanks{Corresponding author: adegoke00@gmail.com}}
\affil{Department of Physics and Engineering Physics, \mbox{Obafemi Awolowo University, Ile-Ife, 220005 Nigeria}}
\author[2]{Olawanle Layeni}
\affil{Department of Mathematics, \mbox{Obafemi Awolowo University, Ile-Ife, 220005 Nigeria}}

\theoremstyle{plain}
\numberwithin{equation}{section}
\newtheorem{thm}{Theorem}[section]

\newtheorem*{theorem}{Theorem}

\newtheorem{example}[thm]{Example}

\newtheorem*{cor}{Corollaries}

\begin{document}

\date{}

\maketitle

\begin{abstract}
\noindent We state and prove a general summation identity. The identity is then applied to derive various summation formulas involving the generalized harmonic numbers and related quantities. Interesting results, mostly new, are obtained for both finite and infinite sums. The high points of this paper are perhaps the discovery of several previously unknown infinite summation results involving {\em non-linear} generalized harmonic number terms and the derivation of interesting alternating summation formulas involving these numbers.
\end{abstract}
\pagebreak
\tableofcontents

\section{Introduction}
Harmonic numbers have been studied since ancient times. Numerous interesting results, especially infinite summation involving these special numbers are scattered in the literature. References~\cite{alzer06, hnmathworld, hnwikipedia, borwein95, flajolet} and further references therein are good sources of information on the subject. In this paper, the generalized harmonic number of order $m$ is denoted by $H_{N,m}$, defined as usual by
\[
H_{N,m}  = \sum\limits_{r = 1}^N {\frac{1}{{r^m }}}\,, 
\]
where $H_{N,1}=H_N$ is the $N-th$ harmonic number. The generalized harmonic number converges to the Riemann Zeta function, $\zeta(m)$:
\[
\lim_{N\to\infty} H_{N,m}=\zeta(m),\quad \Re[m]>1\,,
\]
since
\[
\zeta (m) = \sum\limits_{r = 1}^\infty  {\frac{1}{{r^m }}}\,. 
\]
We define the generalized {\em associated} harmonic number by
\begin{equation}
h_{N,m}  = \sum\limits_{r = 1}^N {\frac{1}{{(2r - 1)^m }}}\,,
\end{equation}
with $h_{N,1}\equiv h_N$ and note that
\[
\lim_{N \to \infty } h_{N,m}  = (1 - 2^{ - m} )\zeta (m),\quad \Re[m]>1\,.
\]
To establish the connection between $H_{N,m}$ and $h_{N,m}$ we first make the following elementary observation:
\begin{equation}\label{equ.n9vef9k}
\sum\limits_{s = 1}^r {f_{s}}  = \sum\limits_{s = 1}^{(r - a_r )/2} {f_{2s}}  + \sum\limits_{s = 1}^{(r + a_r )/2} {f_{2s - 1}}\,,
\end{equation}
where we have introduced the symbol $a_r=r\!\mod 2$.

\bigskip

Taking $f_{s}=1/s^m$ in the identity~\eqref{equ.n9vef9k} allows us to write
\[
H_{r,m}  = \frac{1}{{2^m }}H_{(r - a_r )/2,m}  + \sum\limits_{s = 1}^{(r + a_r )/2} {\frac{1}{{(2s - 1)^m }}}\,, 
\]
which gives, on evaluation at $r=2N$ and at $r=2N-1$, respectively, 
\begin{equation}\label{equ.hroi36f}
\sum\limits_{s = 1}^N {\frac{1}{{(2s - 1)^m }}}  = H_{2N,m}  - \frac{1}{{2^m }}H_{N,m}=h_{N,m}
\end{equation}
and
\begin{equation}\label{equ.n30dcb4}
\sum\limits_{s = 1}^N {\frac{1}{{(2s - 1)^m }}}  = H_{2N - 1,m}  - \frac{1}{{2^m }}H_{N - 1,m}=h_{N,m}\,.
\end{equation}

In what follows, various summation formulas involving $H(r,m)$ and $h(r,m)$ will be derived. Most of these formulas are new and many known results are particular cases of those obtained here. In particular we will derive the following presumably previously unknown summation identities, whose summands contain terms {\em quadratic} in $H(r,2)$, $H(r,3)$ and $h(r,2)$:
\[
\sum\limits_{r = 1}^\infty  {\frac{{H_{r,2}^2 }}{{r^2 }}}  = \frac{{19}}{{22680}}\pi ^6  + \zeta (3)^2,\qquad \sum\limits_{r = 1}^\infty  {\frac{{H_{r,2}^2 }}{{(r+1)^2 }}}  = \frac{{59}}{{22680}}\pi ^6  - \zeta (3)^2\,,
\]
\[
\sum\limits_{r = 1}^\infty  {\frac{{H_{r,2}^2 }}{{r(r + 1)}}}  = \pi ^2\, \zeta (3) - 10\,\zeta (5),\qquad\sum\limits_{r = 1}^\infty  {\frac{{H_{r,3}^2 }}{{r(r + 1)}}}  =  - \frac{{10\pi ^2 }}{3}\,\zeta (5) + 35\,\zeta (7)
\]
and
\[
\sum\limits_{r = 1}^\infty  {\frac{{h_{r,2}^2 }}{{4r^2  - 1}}}  = \frac{{3\pi ^2 }}{{64}}\,\zeta (3)\,.
\]
We will also deduce the following remarkable formulas:
\[
2\sum\limits_{r = 1}^\infty  {( - 1)^{r - 1} H_{r,n} }  = \left( {1 - \frac{1}{{2^{n - 1} }}} \right)\zeta (n),\quad n\ne1\,,
\]
\[
2\sum\limits_{r = 1}^\infty  {( - 1)^{r - 1} h_{r,2n} }  = \beta (2n)\,,\qquad 2\sum\limits_{r = 1}^\infty  {( - 1)^{r - 1} h_{r,2n - 1} }  = \frac{{|E_{2n - 2}| }}{{2^{2n} \Gamma (2n - 1)}}\pi ^{2n - 1}\,,
\]
\[
2\sum\limits_{r = 1}^\infty  {( - 1)^{r - 1} H_{r,n} H_{r - 1,n} }  =  - \frac{{(2^{2n - 1}  - 1)}}{{(2n)!}}\left| {B_{2n} } \right|\pi ^{2n},\qquad 2\sum\limits_{r = 1}^\infty  {( - 1)^{r - 1} h_{r,n} h_{r - 1,n} }  =-\beta {(2n)}\,\,,
\]
where $B_m$ is the $mth$ Bernoulli number, $E_m$ is the $mth$ Euler number and
\[
\beta (m) = \sum\limits_{s = 1}^\infty  {\frac{{( - 1)^{s - 1} }}{{(2s - 1)^m }}}\,. 
\]
Special cases of the above alternating sums include:
\[
2\sum\limits_{r = 1}^\infty  {( - 1)^{r - 1} H_{r,2} }  = \frac{\pi^2}{12},\qquad 2\sum\limits_{r = 1}^\infty  {( - 1)^{r - 1} h_r }  = \frac{\pi }{4}\,,\qquad 2\sum\limits_{r = 1}^\infty  {( - 1)^{r - 1} h_{r,2} }  = G\,,
\]
\[
2\sum\limits_{r = 1}^\infty  {( - 1)^{r - 1} H_r^2 }  = \frac{{\pi ^2 }}{{12}} - \log ^2 2\,,\qquad 2\sum\limits_{r = 1}^\infty  {( - 1)^{r - 1} h_r h_{r - 1} }  =  - G\,, \qquad 2\sum\limits_{r = 1}^\infty  {( - 1)^{r - 1} h_r^2 }  = \frac{{\pi \log 2}}{4}\,,
\]
where $G=\beta(2)$ is Catalan's constant.

\bigskip

In section~\ref{sec.finite} numerous finite summation formulas will be derived. 

\section{Summation Formula}
\begin{theorem}
Given a non-singular summand, $f_{rs}$, $r,s\in\mathbb{Z^+}$, \mbox{$1\le r,s\le N$}, \mbox{$N\in\mathbb{Z^+}$}, the following summation identity holds:
\begin{equation}\label{equ.umdypab}
\sum_{r=1}^{N}\sum_{s=1}^{r}(f_{rs}+f_{sr})=\sum_{r=1}^N f_{rr}+\sum_{r=1}^{N}\sum_{s=1}^{N}f_{sr}\,.
\end{equation}
\end{theorem}
The proof is by mathematical induction on $N$. The theorem is obviously true for $N=1$. Assume that the proposition is true for $N=K\in\mathbb{Z^+}$, so that 
\[
P_K:\quad \sum_{r=1}^{K}\sum_{s=1}^{r}(f_{rs}+f_{sr})=\sum_{r=1}^K f_{rr}+\sum_{r=1}^{K}\sum_{s=1}^{K}f_{sr}\,.
\]
We now show that $P_{K+1}$ is valid whenever $P_{K}$ holds.
\[
P_{K+1}:\quad \sum_{r=1}^{K+1}\sum_{s=1}^{r}(f_{rs}+f_{sr})=\sum_{r=1}^{K+1} f_{rr}+\sum_{r=1}^{K+1}\sum_{s=1}^{K+1}f_{sr}\,.
\]
\begin{proof}
\[
\begin{split}
&\sum\limits_{r = 1}^{K + 1} {\sum\limits_{s = 1}^r {\left\{ {f_{rs} + f_{sr}} \right\}} }\\ 
&= \sum\limits_{r = 1}^K {\sum\limits_{s = 1}^r {\left\{ {f_{rs} + f_{sr}} \right\}} }  + \sum\limits_{s = 1}^{K + 1} {\left\{ {f_{K+1,s} + f_{s,K+1}} \right\}} \\
& = \sum\limits_{r = 1}^K {\sum\limits_{s = 1}^r {\left\{ {f_{rs} + f_{sr}} \right\}} }  + \sum\limits_{r = 1}^{K + 1} {f_{K+1,r}}  + \sum\limits_{s = 1}^{K + 1} {f_{s,K+1}}\\
&\mbox{ We now invoke $P_K$}\\
&= \sum\limits_{r = 1}^K {\sum\limits_{s = 1}^K {f_{sr}} }  + \sum\limits_{r = 1}^K {f_{rr}}  + \sum\limits_{r = 1}^{K + 1} {f_{K+1,r}}  + \sum\limits_{s = 1}^{K + 1} {f_{s,K+1}}\\
&= \sum\limits_{r = 1}^K {\sum\limits_{s = 1}^K {f_{sr}} } + \sum\limits_{r = 1}^{K + 1} {f_{rr}}+\sum\limits_{r = 1}^K {\sum\limits_{s = K + 1}^{K + 1} {f_{sr}} }  + \sum\limits_{r = K + 1}^{K + 1} {\sum\limits_{s = 1}^{K + 1} {f_{sr}} }\\ 
&= \sum\limits_{r = 1}^K {\sum\limits_{s = 1}^{K + 1} {f_{sr}} }  + \sum\limits_{r = K + 1}^{K + 1} {\sum\limits_{s = 1}^{K + 1} {f_{sr}} }  + \sum\limits_{r = 1}^{K + 1} {f_{rr}}\\
&=\sum\limits_{r = 1}^{K + 1} {\sum\limits_{s = 1}^{K + 1} {f_{sr}} }  + \sum\limits_{r = 1}^{K + 1} {f_{rr}}
\end{split}\,.
\]
\end{proof}

\begin{cor}
\begin{enumerate}
\item[1.] If the summand $f_{rs}$ is symmetric in the summation indices $r$ and $s$, that is, if $f_{rs}=f_{sr}$, then

\begin{equation}\label{equ.j7bav13}
2\sum\limits_{r = 1}^N {\sum\limits_{s = 1}^r {f_{rs}} }  = \sum\limits_{r = 1}^N {f_{rr}}  + \sum\limits_{r = 1}^N {\sum\limits_{s = 1}^N {f_{rs}} }\,.
\end{equation}

\item[2.] If $f_{rs}$ is factorable, that is if $f_{rs}=g_{r}h_{s}$, then

\begin{equation}\label{equ.dx732bj}
\sum\limits_{r = 1}^N {\left\{ {g_{r}\sum\limits_{s = 1}^r {h_{s}} } \right\}}  + \sum\limits_{r = 1}^N {\left\{ {h_{r}\sum\limits_{s = 1}^r {g_{s}} } \right\}}  = \sum\limits_{r = 1}^N {g_{r}h_{r}}  + \left( {\sum\limits_{r = 1}^N {g_{r}} } \right)\left( {\sum\limits_{r = 1}^N {h_{r}} } \right)\,.
\end{equation}
In particular, if $f_{rs}=g_{r}g_{s}$, then
\begin{equation}\label{equ.ln46ypv}
2\sum\limits_{r = 1}^N {\left\{ {g_{r}\sum\limits_{s = 1}^r {g_{s}} } \right\}}  = \sum\limits_{r = 1}^N {\left( {g_{r}} \right)^2 }  + \left( {\sum\limits_{r = 1}^N {g_{r}} } \right)^2\,.
\end{equation}

\item[3.] Setting $f_{rs}=g_{s}$ in identity~\eqref{equ.umdypab} gives

\begin{equation}\label{equ.ip14a9i}
\sum\limits_{r = 1}^N {\sum\limits_{s = 1}^r {g_{s}} }  = (N + 1)\sum\limits_{r = 1}^N {g_{r}}  - \sum\limits_{r = 1}^N {rg_{r}}\,.
\end{equation}

\end{enumerate}

\end{cor}

\section{Applications}

\subsection{General finite summation formulas involving the generalized harmonic numbers}\label{sec.finite}

\begin{example}
Choosing $g_{s}=1/s^n$ in identity~\eqref{equ.ip14a9i} gives

\begin{equation}\label{equ.kb1gfmq}
\sum\limits_{r = 1}^N {H_{r,n} }  = (N + 1)H_{N,n}  - H_{N,n - 1}\,,
\end{equation}

while setting $g_{s}=H_{s,n}$ in identity~\eqref{equ.ip14a9i} and using identity~\eqref{equ.kb1gfmq} gives

\begin{equation}\label{equ.predwns}
2\sum\limits_{r = 1}^N {rH_{r,n} }  = N(N + 1)H_{N,n}  + H_{N,n - 1}  - H_{N,n - 2}\,.
\end{equation}

In particular

\begin{equation}\label{equ.wvl7w5t}
\sum\limits_{r = 1}^N {H_{r} }  = (N + 1)H_{N}  - N\,
\end{equation}

and

\begin{equation}\label{equ.if494a6}
\sum\limits_{r = 1}^N {rH_{r} }  =\frac{1}{2} N(N + 1)H_{N}-\frac{1}{4}N(N-1)\,
\end{equation}

Taking $g_s=sH_{s,n}$ in identity~\eqref{equ.ip14a9i} and using identities~\eqref{equ.kb1gfmq} and \eqref{equ.predwns}, we find

\begin{equation}
\begin{split}
\sum\limits_{r = 1}^N {r^2 H_{r,n} }&= \frac{{N(N + 1)(2N + 1)}}{6}H_{N,n}\\
&\qquad- \frac{1}{6}H_{N,n - 1}  + \frac{1}{2}H_{N,n - 2}  - \frac{1}{3}H_{N,n - 3}\,.
\end{split}
\end{equation} 

In particular
\begin{equation}
\sum\limits_{r = 1}^N {r^2 H_r }  = \frac{{N(N + 1)(2N + 1)}}{6}H_N  - \frac{{N(N - 1)(4N + 1)}}{{36}}\,.
\end{equation}

If we set $g_{s}=H_{s,n}/s^n$ in equation~\eqref{equ.ip14a9i} and make use of equation~\eqref{equ.f2v2vwe}, we obtain the identity
\begin{equation}\label{equ.uo7jljo}
\sum\limits_{r = 1}^N {H_{r,n}^2 }  = (N + 1)H_{N,n}^2  + H_{N,2n - 1}  - 2\sum\limits_{r = 1}^N {\frac{{H_{r,n} }}{{r^{n - 1} }}}\,.
\end{equation}

Upon setting $n=1$ in equation~\eqref{equ.uo7jljo} we obtain the interesting result

\begin{equation}\label{equ.hgtftre}
\sum\limits_{r = 1}^N {H_r^2 }  = (N + 1)H_N^2  - (2N + 1)H_N  + 2N\,.
\end{equation}

Identity~\eqref{equ.kb1gfmq} appeared in~\cite{spiess} (Equation~(43)) and is listed in Wikipedia~\cite{hnwikipedia}. The particular cases, identities~\eqref{equ.wvl7w5t} and \eqref{equ.if494a6} are also derived in~\cite{graham94}, (equation~2.36, page~41 and equation~2.57, page~56).

Using identity~\eqref{equ.n9vef9k}  we write

\[
\sum_{r=1}^N H_{r,n}=\sum_{r=1}^{(N-a_N)/2} H_{2r,n}+\sum_{r=1}^{(N+a_N)/2} H_{2r,n}-\frac{1}{2^n}H_{(N+a_N)/2,n}
\]

from which upon using identity \eqref{equ.kb1gfmq}, we get

\begin{equation}
2\sum\limits_{r = 1}^N {H_{2r,n} }  = 2(N + 1)H_{2N,n}  - H_{2N,n - 1}  - h_{N,n}\,.
\end{equation}

\end{example}

\begin{example}\label{ex.cu17le1}
The choice of $f_{rs}=(2r-1)^{-m}s^{-n}$ in the identity~\eqref{equ.dx732bj} leads to

\begin{equation}\label{equ.afc1dnp}
\sum\limits_{r = 1}^N {\frac{{h_{r,m} }}{{r^n }}}  + \sum\limits_{r = 1}^N {\frac{{H_{r - 1,n} }}{{(2r - 1)^m }}}  = h_{N,m} H_{N,n}\,.
\end{equation}

On setting $n=0$ in identity~\eqref{equ.afc1dnp} we obtain

\begin{equation}\label{equ.f7d4yf6}
\sum\limits_{r = 1}^N {h_{r,m}  = \left( {N + \frac{1}{2}} \right)h_{N,m}  - \frac{1}{2}h_{N,m - 1} }\,.
\end{equation}

In particular,

\begin{equation}\label{equ.ceclnnv}
\sum\limits_{r = 1}^N {h_{r}  = \left( {N + \frac{1}{2}} \right)h_{N}  - \frac{N}{2}}\,.
\end{equation}

Using the identities~\eqref{equ.n9vef9k} and \eqref{equ.gfrtfed} gives

\[
2\sum\limits_{r = 1}^N {h_{r,n} }  = 2\sum\limits_{r = 1}^{(N - a_N )/2} {h_{2r,n} }  + 2\sum\limits_{r = 1}^{(N + a_N )/2} {h_{2r,n} }  - h_{N + a_N ,n}  + \bar h_{N + a_N ,n} 
\]
which, together with identity~\eqref{equ.f7d4yf6} then gives

\begin{equation}
4\sum\limits_{r = 1}^N {h_{2r,n} }  = 2(2N + 1)h_{2N,n}  - h_{2N,n - 1}  - \bar h_{2N,n}\,.
\end{equation}

Substituting $g_s=h_{s,m}$ in identity~\eqref{equ.ip14a9i} and using identity~\eqref{equ.f7d4yf6} gives

\begin{equation}\label{equ.bgfrdse}
\sum\limits_{r = 1}^N {rh_{r,m} }  = \left( {\frac{{N(N + 1)}}{2} + \frac{1}{8}} \right)h_{N,m}  - \frac{1}{8}h_{N,m - 2}\,.
\end{equation}

In particular

\begin{equation}\label{equ.fdjl12h}
\sum\limits_{r = 1}^N {rh_{r} }  = \left( {\frac{{N(N + 1)}}{2} + \frac{1}{8}} \right)h_{N}  - \frac{N^2}{8}\,.
\end{equation}

Taking $g_s=h_{s,n}(2s-1)^{-n}$ in identity~\eqref{equ.ip14a9i} and using the result~\eqref{equ.r00mpzh} we find

\begin{equation}\label{equ.rikv1u5}
2\sum\limits_{r = 1}^N {h_{r,n}^2 }  = (2N + 1)h_{N,n}^2  + h_{N,2n - 1}  - 2\sum\limits_{r = 1}^N {\frac{{h_{r,n} }}{{(2r - 1)^{n - 1} }}}\,.
\end{equation}

Now setting $n=1$ in equation~\eqref{equ.rikv1u5} we obtain

\begin{equation}\label{equ.x0k3m6m}
2\sum\limits_{r = 1}^N {h_r^2 }  = (2N + 1)h_N^2  - 2Nh_N  + N\,.
\end{equation}

\end{example}

\begin{example}

The choice $f_{rs}=H_rH_s$ in identity~\eqref{equ.dx732bj} gives

\begin{equation}\label{equ.gy3kccl}
2\sum\limits_{r = 1}^N {\left\{ {H_r \sum\limits_{s = 1}^r {H_s } } \right\}}  = \sum\limits_{r = 1}^N {H_r^2 }  + \left( {\sum\limits_{r = 1}^N {H_r } } \right)^2\,.
\end{equation}

The use of identities~\eqref{equ.wvl7w5t}, \eqref{equ.if494a6} and \eqref{equ.hgtftre} in identity~\eqref{equ.gy3kccl} leads to

\[
\sum\limits_{r = 1}^N {rH_r^2 }  = \frac{N(N+1)}{2}H_N^2  - \frac{{(N^2 -N - 1)}}{2}H_N  + \frac{N(N-3)}{4}\,.
\]

Similarly, the choice $f_{rs}=h_rh_s$ in identity~\eqref{equ.dx732bj} gives

\begin{equation}\label{equ.xfz3k3z}
2\sum\limits_{r = 1}^N {\left\{ {h_r \sum\limits_{s = 1}^r {h_s } } \right\}}  = \sum\limits_{r = 1}^N {h_r^2 }  + \left( {\sum\limits_{r = 1}^N {h_r } } \right)^2\,.
\end{equation}

The use of identities~\eqref{equ.ceclnnv}, \eqref{equ.fdjl12h} and \eqref{equ.x0k3m6m} in identity~\eqref{equ.xfz3k3z} leads to

\[
\sum\limits_{r = 1}^N {rh_r^2 }  = \frac{{(2N + 1)^2 }}{8}h_N^2  - \frac{{(2N + 1)(2N - 1)}}{{16}}h_N  + \frac{{N^2 }}{{16}}\,.
\]

\end{example}

\begin{example}
Let
\[
f_{r,s}=\frac{x^{pr}y^{qs}}{(r+a)^m(s+b)^n}\,.
\]

$f_{rs}$ is factorable, so we apply equation~\eqref{equ.dx732bj}, which gives immediately

\begin{equation}\label{equ.l7bkfc5}
\begin{split}
&\sum\limits_{r = 1}^N {\left\{ {\frac{{x^{pr} }}{{(r+a)^m }}\sum\limits_{s = 1}^r {\frac{{y^{qs} }}{{(s+b)^n }}} } \right\}}  + \sum\limits_{r = 1}^N {\left\{ {\frac{{y^{qr} }}{{(r+b)^n }}\sum\limits_{s = 1}^r {\frac{{x^{ps} }}{{(s+a)^m }}} } \right\}}\\
&= \sum\limits_{r = 1}^N {\frac{{\left( {x^p y^q } \right)^r }}{{(r+a)^m(r+b)^n }}}  + \left( {\sum\limits_{r = 1}^N {\frac{{x^{pr} }}{{(r+a)^m }}} } \right)\left( {\sum\limits_{r = 1}^N {\frac{{y^{qr} }}{{(r+b)^n }}} } \right)\,.
\end{split}
\end{equation}

Various combinations of the parameters $p,q,m,n,a,b$ and the variables $x, y$ may be considered. As an example if we choose $p=0=q$, then we have the interesting result

\begin{equation}\label{equ.ivtc0sc}
\begin{split}
&\sum\limits_{r = 1}^N {\frac{{H_{r + b,n} }}{{(r + a)^m }}}+\sum\limits_{r = 1}^N {\frac{{H_{r + a,m} }}{{(r + b)^n }}}\\
&\qquad = H_{N + a,m} H_{N + b,n}  - H_{a,m} H_{b,n}  + \sum\limits_{r = 1}^N {\frac{1}{{(r + a)^m (r + b)^n }}}\,.
\end{split}
\end{equation}

In deriving the identity~\eqref{equ.ivtc0sc} we made use of the identity

\[
\sum\limits_{t = 1}^s {\frac{1}{{(t + q)^p }}}  = H_{s + q,p}  - H_{q,p}\,. 
\]

Interesting special cases of identity~\eqref{equ.ivtc0sc} include

\begin{equation}\label{equ.ybdhkk3}
\sum\limits_{r = 1}^N {\frac{{H_{r,n} }}{{r^m }}}  + \sum\limits_{r = 1}^N {\frac{{H_{r,m} }}{{r^n }}}  = H_{N,m + n}  + H_{N,m} H_{N,n}\,,
\end{equation}

\begin{equation}\label{equ.itn3i7s}
\begin{split}
&\sum\limits_{r = 1}^N {\frac{{H_{r,n} }}{{(r + 1)^m }}}  + \sum\limits_{r = 1}^N {\frac{{H_{r,m} }}{{(r + 1)^n }}}\\
&\quad = H_{N + 1,n} H_{N + 1,m}  - H_{N,n + m}  - \frac{1}{{(N + 1)^n (N + 1)^m }}
\end{split}
\end{equation}

and

\begin{equation}
\sum\limits_{r = 1}^N {\frac{{H_{r,n} }}{{(r + 1)^m }}}  + \sum\limits_{r = 1}^N {\frac{{H_{r,m} }}{{r^n }}}  = H_{N + 1,m} H_{N,n}\,.
\end{equation}

The particular case $m=n$ in equations~\eqref{equ.ybdhkk3} and \eqref{equ.itn3i7s} gives

\begin{equation}\label{equ.f2v2vwe}
2\sum\limits_{r = 1}^N {\frac{{H_{r,n} }}{{r^n }}}  = H_{N,2n}  + H_{N,n}^2
\end{equation}
and
\begin{equation}\label{equ.quhf1ku}
2\sum\limits_{r = 1}^N {\frac{{H_{r,n} }}{{(r + 1)^n }}}  = H_{N,n}^2  - H_{N,2n}  + \frac{{2H_{N,n} }}{{(N + 1)^n }}\,.
\end{equation}

The particular case corresponding to $n=1$ in~\eqref{equ.quhf1ku} is also found in~\cite{spiess} (page 850, Theorem~16, example).

\bigskip

Equation (3.62) of reference~\cite{alzer06} corresponds to setting $n=1$ in identity~\eqref{equ.f2v2vwe}. 

\end{example}

\begin{example}\label{ex.e363h43}

Substitution of $f_{rs}=(2r+2a-1)^{-m}(2s+2b-1)^{-n}$ into equation~\eqref{equ.dx732bj} gives

\begin{equation}\label{equ.d7w4dco}
\begin{split}
&\sum\limits_{r = 1}^N {\frac{{h_{r + b,n} }}{{(2r + 2a-1)^m }}}+\sum\limits_{r = 1}^N {\frac{{h_{r + a,m} }}{{(2r + 2b-1)^n }}}\\
&\qquad = h_{N + a,m} h_{N + b,n}  - h_{a,m} h_{b,n}  + \sum\limits_{r = 1}^N {\frac{1}{{(2r + 2a-1)^m (2r + 2b-1)^n }}}\,.
\end{split}
\end{equation}

Note that in deriving the identity~\eqref{equ.d7w4dco} we made use of the identity

\[
\sum\limits_{t = 1}^s {\frac{1}{{(2t + 2q-1)^p }}}  = h_{s + q,p}  - h_{q,p}\,. 
\]

Interesting special cases of identity~\eqref{equ.d7w4dco} include

\begin{equation}\label{equ.vse99ju}
\sum\limits_{r = 1}^N {\frac{{h_{r,n} }}{{(2r - 1)^m }}}  + \sum\limits_{r = 1}^N {\frac{{h_{r,m} }}{{(2r - 1)^n }}}  = h_{N,n + m}  + h_{N,m} h_{N,n}\,,
\end{equation}

\begin{equation}\label{equ.nrl5ek4}
\begin{split}
&\sum\limits_{r = 1}^N {\frac{{h_{r,n} }}{{(2r + 1)^m }}}  + \sum\limits_{r = 1}^N {\frac{{h_{r,m} }}{{(2r + 1)^n }}}\\
&\quad = h_{N + 1,n} h_{N + 1,m}  - h_{N,n + m}  - \frac{1}{{(2N + 1)^n (2N + 1)^m }}
\end{split}
\end{equation}

and

\begin{equation}\label{equ.hyn439v}
\sum\limits_{r = 1}^N {\frac{{h_{r,n} }}{{(2r + 1)^m }}}  + \sum\limits_{r = 1}^N {\frac{{h_{r,m} }}{{(2r-1)^n }}}  = h_{N + 1,m} h_{N,n}\,.
\end{equation}

The particular case $m=n$ in equations~\eqref{equ.vse99ju} and \eqref{equ.nrl5ek4} gives

\begin{equation}\label{equ.r00mpzh}
2\sum\limits_{r = 1}^N {\frac{{h_{r,n} }}{{(2r-1)^n }}}  = h_{N,2n}  + h_{N,n}^2
\end{equation}

and

\begin{equation}\label{equ.wfa9yp1}
2\sum\limits_{r = 1}^N {\frac{{h_{r,n} }}{{(2r + 1)^n }}}  = h_{N,n}^2  - h_{N,2n}  + \frac{{2h_{N,n} }}{{(2N + 1)^n }}\,.
\end{equation}

From identities~\eqref{equ.r00mpzh} and \eqref{equ.wfa9yp1} we have

\begin{equation}
\sum\limits_{r = 1}^N {\frac{{h_{r,n} }}{{(2r-1)^n }}}  + \sum\limits_{r = 1}^N {\frac{{h_{r,n} }}{{(2r + 1)^n }}}  = h_{N,n}^2  + \frac{{h_{N,n} }}{{(2N + 1)^n }}
\end{equation}
and
\begin{equation}
\sum\limits_{r = 1}^N {\frac{{h_{r,n} }}{{(2r-1)^n }}}  - \sum\limits_{r = 1}^N {\frac{{h_{r,n} }}{{(2r + 1)^n }}}  = h_{N,2n}  - \frac{{h_{N,n} }}{{(2N + 1)^n }}\,.
\end{equation}

Again all the formulas derived in this example are new.

\end{example}

\begin{example}\label{ex.bmlp79d}

Substitution of $f_{rs}=r^{-n}z^{s}$ into identity~\eqref{equ.dx732bj} gives, after some rearrangement,

\begin{equation}\label{equ.bnils61}
\sum\limits_{r = 1}^N {z^r H_{r,n} }  = \frac{1}{{1 - z}}\sum\limits_{r = 1}^N {\frac{{z^r }}{{r^n }}}  - \frac{{z^{N + 1} }}{{1 - z}}H_{N,n},\quad z\ne 1\,,
\end{equation}

while substitution of $f_{rs}=(2r-1)^{-n}z^{(2s-1)}$ into identity~\eqref{equ.dx732bj} yields

\begin{equation}\label{equ.p5u5zed}
\sum\limits_{r = 1}^N {z^{2r - 1} h_{r,n} }  = \frac{1}{{1 - z^2 }}\sum\limits_{r = 1}^N {\frac{{z^{2r - 1} }}{{(2r - 1)^n }}}  - \frac{{z^{2N + 1} }}{{1 - z^2 }}h_{N,n},\quad z\ne 1\,.
\end{equation}
\end{example}

\begin{example}

If we choose $f_{rs}=H_{r,n}/{r^ns^m}$ in the equation~\eqref{equ.dx732bj} we obtain the following identity, valid for all complex numbers $n,m$ and positive integers $N$:
\begin{equation}
\begin{split}
&2\sum\limits_{r = 1}^N {\frac{{H_{r,n} H_{r,m} }}{{r^n }}}  + \sum\limits_{r = 1}^N {\frac{{H_{r,2n} }}{{r^m }}}  + \sum\limits_{r = 1}^N {\frac{{H_{r,n}^2 }}{{r^m }}}\\
& \qquad\qquad= 2\sum\limits_{r = 1}^N {\frac{{H_{r,n} }}{{r^{m + n} }}}  + H_{N,m} H_{N,2n}  + H_{N,m} H_{N,n}^2\,.
\end{split}
\end{equation}

In particular, setting $m=n$ and using also the identity~\eqref{equ.f2v2vwe} we obtain the beautiful result

\begin{equation}\label{equ.ubtcdej}
3\sum\limits_{r = 1}^N {\frac{{H_{r,n}^2 }}{{r^n }}}  - 3\sum\limits_{r = 1}^N {\frac{{H_{r,n} }}{{r^{2n} }}}  = H_{N,n}^3  - H_{N,3n}\,,
\end{equation}

or equivalently,

\begin{equation}\label{equ.bo6qz7a}
3\sum\limits_{r = 1}^N {\frac{{H_{r,n}^2 }}{{r^n }}}  + 3\sum\limits_{r = 1}^N {\frac{{H_{r,2n} }}{{r^{n} }}}  = H_{N,n}^3+3H_{N,2n}H_{N,n}  +2 H_{N,3n}\,.
\end{equation}

Note that since

\[
H_{r,n}  = H_{r + 1,n}  - \frac{1}{{(r + 1)^n }}\,,
\]

identity~\eqref{equ.ubtcdej} can also be written

\begin{equation}\label{equ.mimii8d}
3\sum\limits_{r = 1}^N {\frac{{H_{r,n}^2 }}{{(r + 1)^n }}}  + 3\sum\limits_{r = 1}^N {\frac{{H_{r,n} }}{{r^{2n} }}}  = H_{N + 1,n}^3  + 2H_{N,3n}  + \frac{2}{{(N + 1)^{3n} }} - \frac{{3H_{N + 1,n} }}{{(N + 1)^{2n} }}\,.
\end{equation}

Addition of identities~\eqref{equ.ubtcdej} and \eqref{equ.mimii8d} gives

\begin{equation}\label{equ.h42rpes}
\begin{split}
3\sum\limits_{r = 1}^N {\frac{{H_{r,n}^2 }}{{(r + 1)^n }}}  + 3\sum\limits_{r = 1}^N {\frac{{H_{r,n}^2 }}{{r^n }}}  &= H_{N + 1,n}^3  + H_{N,n}^3  + H_{N,3n}\\
&\qquad + \frac{2}{{(N + 1)^{3n} }} - \frac{{3H_{N + 1,n} }}{{(N + 1)^{2n} }}\,.
\end{split}
\end{equation}

\end{example}

\begin{example}
The choice $f_{rs}=h_{r,n}{(2r-1)^{-n}(2s-1)^{-m}}$ in equation~\eqref{equ.dx732bj} yields the following identity, which holds for all complex numbers $n,m$ and positive integers $N$:
\begin{equation}
\begin{split}
&2\sum\limits_{r = 1}^N {\frac{{h_{r,n} h_{r,m} }}{{(2r-1)^n }}}  + \sum\limits_{r = 1}^N {\frac{{h_{r,2n} }}{{(2r-1)^m }}}  + \sum\limits_{r = 1}^N {\frac{{h_{r,n}^2 }}{{(2r-1)^m }}}\\
& \qquad\qquad= 2\sum\limits_{r = 1}^N {\frac{{h_{r,n} }}{{(2r-1)^{m + n} }}}  + h_{N,m} h_{N,2n}  + h_{N,m} h_{N,n}^2\,.
\end{split}
\end{equation}

In particular, setting $m=n$ and using also the identity~\eqref{equ.r00mpzh} we obtain the interesting result

\begin{equation}\label{equ.mr00dwg}
3\sum\limits_{r = 1}^N {\frac{{h_{r,n}^2 }}{{(2r-1)^n }}}  - 3\sum\limits_{r = 1}^N {\frac{{h_{r,n} }}{{(2r-1)^{2n} }}}  = h_{N,n}^3  - h_{N,3n}\,,
\end{equation}

or equivalently,

\begin{equation}\label{equ.jkfa132}
3\sum\limits_{r = 1}^N {\frac{{h_{r,n}^2 }}{{(2r-1)^n }}}  + 3\sum\limits_{r = 1}^N {\frac{{h_{r,2n} }}{{(2r-1)^{n} }}}  = h_{N,n}^3+3h_{N,2n}h_{N,n}  +2 h_{N,3n}\,.
\end{equation}

\end{example}

\begin{example}
In this example we derive a couple of alternating summation formulas.

\bigskip
First we introduce the notations
\begin{equation}\label{equ.ol6x7b4}
\bar H_{p,q}  = \sum\limits_{s = 1}^p {\frac{{( - 1)^{s - 1} }}{{s^q }}}\text{ and } \bar h_{p,q}  = \sum\limits_{s = 1}^p {\frac{{( - 1)^{s - 1} }}{{(2s - 1)^q }}}\,.
\end{equation}

Then, from the identity~\eqref{equ.n9vef9k} we have

\[
\bar H_{N,n}  =  - \frac{1}{{2^n }}H_{(N - a_N )/2,n}  + h_{(N + a_N )/2,n}\,, 
\]

from which it follows that

\begin{equation}\label{equ.to49yt7}
\bar H_{2N,n}  =  - \frac{1}{{2^n }}H_{N,n}  + h_{N,n}
\end{equation}

and

\begin{equation}\label{equ.ldfid07}
\bar H_{2N - 1,n}  =  - \frac{1}{{2^n }}H_{N - 1,n}  + h_{N,n}\,.
\end{equation}

Similarly, using the identity~\eqref{equ.n9vef9k} and the definitions of $h$ and $\bar h$, it is straightforward to establish that

\begin{equation}\label{equ.gfrtfed}
2\sum\limits_{r = 1}^N {\frac{1}{{(4r - 1)^n }}}  = h_{2N,n}  - \bar h_{2N,n}= \frac{{2( - 1)^{n - 1} }}{{4^n \Gamma (n)}}\left\{ {\psi _{n - 1} \left( {N + \frac{3}{4}} \right) - \psi _{n - 1} \left( {\frac{3}{4}} \right)} \right\}
\end{equation}

and

\begin{equation}\label{equ.nhgftfj}
2\sum\limits_{r = 1}^N {\frac{1}{{(4r - 3)^n }}}  = h_{2N,n}  + \bar h_{2N,n} = \frac{{2( - 1)^{n - 1} }}{{4^n \Gamma (n)}}\left\{ {\psi _{n - 1} \left( {N + \frac{1}{4}} \right) - \psi _{n - 1} \left( {\frac{1}{4}} \right)} \right\}\,,
\end{equation}

where $\psi_n(x)$ is the $nth$ polygamma function defined by \[ \psi_n(x)=\frac{d\psi(x)}{dx^n}\] where \[\psi(x)=\frac{d}{dx}\log\Gamma(x)\] is the digamma function and $\Gamma(x)$ is the gamma function.

Using $f_{rs}=(-1)^{s-1}r^{-n}$ in identity~\eqref{equ.dx732bj} we obtain

\begin{equation}
\sum_{r=1}^N (-1)^{r-1} H_{r,n}=-\frac{1}{2^n}H_{\frac{N-a_N}{2},n}+a_NH_{N,n}
\end{equation}

from which we get the interesting results

\begin{equation}\label{equ.h87qnmg}
\sum_{r=1}^{2N} (-1)^{r-1} H_{r,n}=-\frac{1}{2^n}H_{N,n}
\end{equation}

and

\begin{equation}\label{equ.bx3kgp2}
\sum_{r=1}^{2N-1} (-1)^{r-1} H_{r,n}=h_{N,n}\,.
\end{equation}

Similarly using $f_{rs}=(-1)^{s-1}(2r-1)^{-n}$ in identity~\eqref{equ.dx732bj} gives

\[
\sum\limits_{r = 1}^N {( - 1)^{r - 1} h_{r,n} }  =  - \sum\limits_{r = 1}^{(N - a_N )/2} {\frac{1}{{(4r - 1)^n }}}  + a_N h_{N,n}\,, 
\]
which leads to

\begin{equation}\label{equ.n7gklxg}
2\sum\limits_{r = 1}^{2N} {( - 1)^{r - 1} h_{r,n} }  =  \bar h_{2N,n}-h_{2N,n}
\end{equation}
and

\begin{equation}\label{equ.idbeoly}
2\sum\limits_{r = 1}^{2N - 1} {( - 1)^{r - 1} h_{r,n} }  = \bar h_{2N,n}+h_{2N,n}\,.
\end{equation}

The particular case corresponding to $n=1$ in identity~\eqref{equ.h87qnmg} is also derived in~\cite{spiess} (Equation~(39)).

\bigskip

Using $f_{rs}=(-1)^{r-1}(-1)^{s-1}r^{-n}$ in identity~\eqref{equ.dx732bj} yields

\begin{equation}\label{equ.kv6sa2t}
\sum\limits_{r = 1}^{2N} {( - 1)^{r - 1} \bar H_{r,n} }  = {{H_{N,n} } \mathord{\left/
 {\vphantom {{H_{N,n} } {2^n }}} \right.
 \kern-\nulldelimiterspace} {2^n }}
\end{equation}

and

\begin{equation}\label{equ.um6wh8z}
\sum\limits_{r = 1}^{2N - 1} {( - 1)^{r - 1} \bar H_{r,n} }  = h_{N,n}\,.
\end{equation}

Taking $f_{rs}=(-1)^{(s-1)}H_{s,n}r^{-m}$ in identity~\eqref{equ.dx732bj} gives

\begin{equation}\label{equ.an6o3qx}
\begin{split}
&- \frac{1}{{2^{m + n} }}\sum\limits_{r = 1}^{(N - a_N )/2} {\frac{{H_{r,n} }}{{r^m }}}  + \sum\limits_{r = 1}^{(N + a_N )/2} {\frac{{h_{r,n} }}{{(2r - 1)^m }}}\\
&\qquad+ \sum\limits_{r = 1}^N {( - 1)^{r - 1} H_{r,n} H_{r,m} }\\
&\qquad\qquad= \sum\limits_{r = 1}^N {( - 1)^{r - 1} \frac{{H_{r,n} }}{{r^m }}}  + H_{N,m} \sum\limits_{r = 1}^N {( - 1)^{r - 1} H_{r,n} }\,.
\end{split}
\end{equation}

Interchanging $m$ and $n$ in identity~\eqref{equ.an6o3qx}, adding the resulting identity to identity~\eqref{equ.an6o3qx} and using identities~\eqref{equ.ybdhkk3} and \eqref{equ.vse99ju} we obtain

\begin{equation}
\begin{split}
& - \frac{1}{{2^{m + n} }}\left( {H_{(N - a_N )/2,n + m}  + H_{(N - a_N )/2,n} H_{(N - a_N )/2,m} } \right)\\
&\quad + h_{(N + a_N )/2,n + m}  + h_{(N + a_N )/2,n} h_{(N + a_N )/2,m}\\
&\qquad + 2\sum\limits_{r = 1}^N {( - 1)^{r - 1} H_{r,n} H_{r,m} }\\
&\quad\qquad = \sum\limits_{r = 1}^N {( - 1)^{r - 1} \left( {\frac{{H_{r,n} }}{{r^m }} + \frac{{H_{r,m} }}{{r^n }}} \right)}\\
&\qquad\qquad + H_{N,m} \sum\limits_{r = 1}^N {( - 1)^{r - 1} H_{r,n} }\\
&\qquad\qquad + H_{N,n} \sum\limits_{r = 1}^N {( - 1)^{r - 1} H_{r,m} }\,,
\end{split}
\end{equation}

from which we finally get

\begin{equation}\label{equ.n8tkvqw}
\begin{split}
&\sum\limits_{r = 1}^{2N} {( - 1)^{r - 1} \left( {2H_{r,n} H_{r,m}  - \frac{{H_{r,n} }}{{r^m }} - \frac{{H_{r,m} }}{{r^n }}} \right)}\\
&\quad = \frac{1}{{2^{m + n} }}\left( {H_{N,m + n}  + H_{N,m} H_{N,n} } \right)\\
&\qquad- h_{N,m + n}  - h_{N,m} h_{N,n}\\
&\quad\qquad - \frac{{H_{2N,m} H_{N,n} }}{{2^n }} - \frac{{H_{2N,n} H_{N,m} }}{{2^m }}
\end{split}
\end{equation}

and

\begin{equation}\label{equ.av4l176}
\begin{split}
&\sum\limits_{r = 1}^{2N-1} {( - 1)^{r - 1} \left( {2H_{r,n} H_{r,m}  - \frac{{H_{r,n} }}{{r^m }} - \frac{{H_{r,m} }}{{r^n }}} \right)}\\
&\quad = \frac{1}{{2^{m + n} }}\left( {H_{N-1,m + n}  + H_{N-1,m} H_{N-1,n} } \right)\\
&\qquad- h_{N,m + n}  - h_{N,m} h_{N,n}\\
&\quad\qquad +H_{2N-1,m} h_{N,n}+H_{2N-1,n} h_{N,m}\,.
\end{split}
\end{equation}

In particular

\begin{equation}
2\sum\limits_{r = 1}^{2N} {( - 1)^{r - 1} rH_{r,n} }  = 2h_{N,n - 1}  - h_{N,n}  - 2NH_{2N,n}  - H_{2N,n - 1}\,,
\end{equation}

\begin{equation}
2\sum\limits_{r = 1}^{2N - 1} {( - 1)^{r - 1} rH_{r,n} }  = 2h_{N,n - 1}  - h_{N,n}  + 2NH_{2N - 1,n}  - H_{2N - 1,n - 1}\,,
\end{equation}

\begin{equation}\label{equ.mk2d5hz}
\sum\limits_{r = 1}^{2N} {( - 1)^{r - 1} \left( {2H_{r,n}^2  - 2\frac{{H_{r,n} }}{{r^n }}} \right)}  = \frac{{H_{N,2n} }}{{2^{2n} }} - h_{N,2n}  - H_{2N,n}^2
\end{equation}

and

\begin{equation}\label{equ.gn5ox3z}
\sum\limits_{r = 1}^{2N - 1} {( - 1)^{r - 1} \left( {2H_{r,n}^2  - 2\frac{{H_{r,n} }}{{r^n }}} \right)}  = \frac{{H_{N - 1,2n} }}{{2^{2n} }} - h_{N,2n}  + H_{2N - 1,n}^2\,.
\end{equation}

Corresponding to identities~\eqref{equ.mk2d5hz} and \eqref{equ.gn5ox3z} we have, upon taking \mbox{$f_{rs}=(-1)^{(s-1)}h_{s,n}(2r-1)^{-n}$} in identity~\eqref{equ.dx732bj}

\begin{equation}\label{equ.gpqd08z}
2\sum\limits_{r = 1}^{2N} {( - 1)^{r - 1} \left( {h_{r,n}^2  - \frac{{h_{r,n} }}{{(2r - 1)^n }}} \right)}  =  - h_{2N,n}^2  - \bar h_{2N,2n}
\end{equation}

and

\begin{equation}\label{equ.v1b4oyo}
2\sum\limits_{r = 1}^{2N - 1} {( - 1)^{r - 1} \left( {h_{r,n}^2  - \frac{{h_{r,n} }}{{(2r - 1)^n }}} \right)}  = h_{2N - 1,n}^2  - \bar h_{2N - 1,2n}\,.
\end{equation}

\end{example}

\subsection{Evaluation of infinite sums}

In the limit $N\to\infty$ in the above summation results and sometimes in combination with known results, it is possible to evaluate certain infinite sums. We now present some examples.

\begin{example}

In the limit $N\to\infty$, equations~\eqref{equ.ubtcdej}, \eqref{equ.bo6qz7a} and \eqref{equ.htujppt} become

\begin{equation}\label{equ.jx9a49p}
3\sum\limits_{r = 1}^\infty  {\frac{{H_{r,n}^2 }}{{r^n }}}  - 3\sum\limits_{r = 1}^\infty  {\frac{{H_{r,n} }}{{r^{2n} }}}  = \zeta (n)^3  - \zeta (3n),\quad n\ne 1\,,
\end{equation}

\begin{equation}
3\sum\limits_{r = 1}^\infty  {\frac{{H_{r,n}^2 }}{{r^n }}}  + 3\sum\limits_{r = 1}^\infty  {\frac{{H_{r,2n} }}{{r^n }}}  = \zeta (n)^3  + 3\zeta (n)\zeta (2n) + 2\zeta (3n),\quad n\ne 1
\end{equation}

and

\begin{equation}\label{equ.b7z8vy9}
3\sum\limits_{r = 1}^\infty  {\frac{{H_{r,n}^2 }}{{r^n }}} + 3\sum\limits_{r = 1}^\infty  {\frac{{H_{r,n}^2 }}{{(r+1)^n }}}  = 2\zeta (n)^3  + \zeta (3n),\quad n\ne 1\,.
\end{equation}

Evaluating identity~\eqref{equ.jx9a49p} at $n=2$ we obtain

\begin{equation}\label{equ.htujppt}
\sum\limits_{r = 1}^\infty  {\frac{{H_{r,2}^2 }}{{r^2 }}}  = \frac{{19}}{{22680}}\pi ^6  + \zeta (3)^2\,,
\end{equation}
after using the known result:
\[
\begin{split}
\sum\limits_{r = 1}^\infty  {\frac{{H_{r,2} }}{{r^4 }}}  = \zeta (3)^2  - \frac{{\pi ^6 }}{{2835}},\quad \mbox{(\cite{coffey}, (B.9a), \cite{hnmathworld})}\,.
\end{split}
\]

Now using the result~\eqref{equ.htujppt} in identity~\eqref{equ.b7z8vy9}, we also have

\begin{equation}\label{equ.xx2cqv4}
\sum\limits_{r = 1}^\infty  {\frac{{H_{r,2}^2 }}{{(r+1)^2 }}}  = \frac{{59}}{{22680}}\pi ^6  - \zeta (3)^2\,.
\end{equation}

Since
\[
H_{r - 1,n}^2  = \left( {H_{r,n}  - \frac{1}{{r^n }}} \right)^2  = H_{r,n}^2  - \frac{{2H_{r,n} }}{{r^n }} + \frac{1}{{r^{2n} }}\,,
\]
we have

\[
\sum\limits_{r = 1}^\infty  {\frac{{H_{r,2}^2 }}{{r(r + 1)}}}  = 2\sum\limits_{r = 1}^\infty  {\frac{{H_{r,2} }}{{r^3 }}}  - \zeta (5)\,,
\]
and
\[
\sum\limits_{r = 1}^\infty  {\frac{{H_{r,3}^2 }}{{r(r + 1)}}}  = 2\sum\limits_{r = 1}^\infty  {\frac{{H_{r,3} }}{{r^4 }}}  - \zeta (7)\,,
\]

from which upon using the known results

\[
2\sum\limits_{r = 1}^\infty  {\frac{{H_{r,2} }}{{r^3 }}}  = \pi ^2 \zeta (3) - 9\zeta (5),\quad \mbox{(Eq. 3.3b of \cite{zheng})}\,,
\]
and
\[
\sum\limits_{r = 1}^\infty  {\frac{{H_{r,4} }}{{r^3 }}}  = \frac{{\pi ^4 }}{{90}}\zeta (3) - \frac{{5\pi ^2 }}{3}\zeta (5) - 17\zeta (7),\quad \mbox{(Eq. 3.5d of \cite{zheng})}\,,
\]

we obtain

\[
\sum\limits_{r = 1}^\infty  {\frac{{H_{r,2}^2 }}{{r(r + 1)}}}  = \pi ^2 \zeta (3) - 10\zeta (5)
\]

and

\[
\sum\limits_{r = 1}^\infty  {\frac{{H_{r,3}^2 }}{{r(r + 1)}}}  = 35\zeta (7) - \frac{10\pi^2}{3}\zeta (5)\,.
\]

\end{example}

\begin{example}

In the limit $N\to\infty$, equations~\eqref{equ.mr00dwg} and \eqref{equ.jkfa132} become

\begin{equation}\label{equ.wc5idtu}
3\sum\limits_{r = 1}^\infty {\frac{{h_{r,n}^2 }}{{(2r-1)^n }}}  - 3\sum\limits_{r = 1}^\infty {\frac{{h_{r,n} }}{{(2r-1)^{2n} }}}  = (1-2^{-n})^3\zeta(n)^3  - (1-2^{-3n})\zeta({3n})\,,
\end{equation}
and
\begin{equation}\label{equ.a710hii}
\begin{split}
&3\sum\limits_{r = 1}^\infty {\frac{{h_{r,n}^2 }}{{(2r-1)^n }}}  + 3\sum\limits_{r = 1}^\infty {\frac{{h_{r,2n} }}{{(2r-1)^{n} }}}\\
&\quad  = (1-2^{-n})^3\zeta(n)^3 +3(1-2^{-n})(1-2^{-2n})\zeta(n)\zeta(2n)  +2 (1-2^{-3n})\zeta(3n)\,.
\end{split}
\end{equation}

\end{example}

\begin{example}

Dividing through identity~\eqref{equ.f2v2vwe} by $r^m$, summing and taking limit as $N\to\infty$ gives

\begin{equation}\label{equ.kj6can4}
2\sum\limits_{r = 1}^\infty  {\left\{ {\frac{1}{{r^m }}\sum\limits_{s = 1}^r {\frac{{H_{s,n} }}{{s^n }}} } \right\}}  = \sum\limits_{r = 1}^\infty  {\frac{{H_{r,2n} }}{{r^m }}}  + \sum\limits_{r = 1}^\infty  {\frac{{H_{r,n}^2 }}{{r^m }}},\quad m\ne 1\,.
\end{equation}

In particular $(m,n)=(2,1)$ and $(m,n)=(2,2)$ in~\eqref{equ.kj6can4} give, respectively,

\begin{equation}\label{equ.tz92kxc}
2\sum\limits_{r = 1}^\infty  {\left\{ {\frac{1}{{r^2 }}\sum\limits_{s = 1}^r {\frac{{H_s }}{s}} } \right\}}  = \sum\limits_{r = 1}^\infty  {\frac{{H_{r,2} }}{{r^2 }}}  + \sum\limits_{r = 1}^\infty  {\frac{{H_r^2 }}{{r^2 }}}
\end{equation}

and

\begin{equation}\label{equ.uyj3ssj}
2\sum\limits_{r = 1}^\infty  {\left\{ {\frac{1}{{r^2 }}\sum\limits_{s = 1}^r {\frac{{H_{s,2} }}{{s^2 }}} } \right\}}  = \sum\limits_{r = 1}^\infty  {\frac{{H_{r,4} }}{{r^2 }}}  + \sum\limits_{r = 1}^\infty  {\frac{{H_{r,2}^2 }}{{r^2 }}}\,.
\end{equation}

Using equation~\eqref{equ.khvarp3} evaluated at $n=2$ and the known result
\[
\sum\limits_{r = 1}^\infty  {\frac{{H_r^2 }}{{r^2 }}}  = \frac{{17}}{{360}}\pi ^4,\quad\mbox{(\cite{borwein95}, \cite{alzer06})}\,, 
\]

in equation~\eqref{equ.tz92kxc} we obtain

\begin{equation}\label{equ.rhge266}
\sum\limits_{r = 1}^\infty  {\left\{ {\frac{1}{{r^2 }}\sum\limits_{s = 1}^r {\frac{{H_s }}{s}} } \right\}}=\frac{\pi^4}{30}\,.
\end{equation}

Using the result~\eqref{equ.htujppt} above and the known result

\[
\sum\limits_{r = 1}^\infty  {\frac{{H_{r,4} }}{{r^2 }}}  = \frac{{37}}{{11340}}\pi ^6  - \zeta (3)^2,\quad\mbox{(Formula~(42) of~\cite{hnmathworld})}\,, 
\]
in equation~\eqref{equ.uyj3ssj} we obtain

\begin{equation}\label{equ.e6pzke8}
\sum\limits_{r = 1}^\infty  {\left\{ {\frac{1}{{r^2 }}\sum\limits_{s = 1}^r {\frac{{H_{s,2} }}{{s^2 }}} } \right\}}  = \frac{{31}}{{15120}}\pi ^6\,.
\end{equation}

Equation~\eqref{equ.rhge266} was also derived in reference~\cite{alzer06}.

\end{example}

\begin{example}

In the limit $N\to\infty$ in equation~\eqref{equ.bnils61} of Example~\ref{ex.bmlp79d}, we get the known result (Formula~(36) of~\cite{hnmathworld})

\begin{equation}\label{equ.bzkcrl5}
\sum\limits_{r = 1}^\infty  z^rH_{r,n} = \frac{1}{{1 - z}}\Li_n(z)\,, |z|<1\,,
\end{equation}

where $\Li_n$ is the polylogarithm function.

\bigskip

At $n=1$ we have

\[
\sum\limits_{r = 1}^\infty  {z^r H_r }  =  - \frac{{\log (1 - z)}}{{1 - z}},\quad |z|<1\,.
\]

Other interesting particular cases are

\[
\sum\limits_{r = 1}^\infty  {\frac{{H_{r,2} }}{{2^r }}}  = \frac{{\pi ^2 }}{6} - \log ^2 2
\]

and

\[
\sum\limits_{r = 1}^\infty  {\frac{{H_{r,3} }}{{2^r }}}  = \frac{7}{4}\zeta (3) - \frac{1}{6}\pi ^2 \log ^2 2 + \frac{1}{3}\log ^3 2\,.
\]

\bigskip

Using the recurrence relation of the polygamma function

\begin{equation}\label{equ.m2dzb1g}
\psi _m (z + 1) = \psi _m (z) + \frac{{( - 1)^m m!}}{{z^{m + 1} }}
\end{equation}

and the identity

\begin{equation}\label{equ.ybh1wcb}
\frac{{\psi _m (z)}}{{( - 1)^{m + 1} m!}} = \zeta (m + 1) - H_{z - 1,m + 1}\,,
\end{equation}

equation~\eqref{equ.bzkcrl5} can be written in terms of the polygamma function as

\[
\sum\limits_{r = 1}^\infty  {z^r \psi _{n - 1} (r)}  = ( - 1)^{(n-1)} (n - 1)!\frac{z}{{1 - z}}\left[ {\Li_n (z) - \zeta (n)} \right],\quad n>1,\,|z|<1\,.
\]

In the limit $N\to\infty$, identity~\eqref{equ.p5u5zed} becomes

\[
2\sum\limits_{r = 1}^\infty  {z^{2r - 1} h_{r,n} }  = \frac{{\Li_n (z) - \Li_n ( - z)}}{{1 - z^2 }},\quad |z|<1\,.
\]

In particular,

\[
2\sum\limits_{r = 1}^\infty  {z^{2r - 1} h_r }  = \frac{1}{{1 - z^2 }}\log \left( {\frac{{1 + z}}{{1 - z}}} \right),\quad |z|<1\,.
\]

\end{example}

\begin{example}

In the limit of $N\to\infty$, equation~\eqref{equ.l7bkfc5} becomes

\begin{equation}\label{equ.o6mnea9}
\begin{split}
\sum\limits_{r = 1}^\infty {\left\{ {\frac{{x^{pr} }}{{r^m }}\sum\limits_{s = 1}^r {\frac{{y^{qs} }}{{s^n }}} } \right\}}  + \sum\limits_{r = 1}^\infty {\left\{ {\frac{{y^{qr} }}{{r^n }}\sum\limits_{s = 1}^r {\frac{{x^{ps} }}{{s^m }}} } \right\}}=\Li_{m + n} \left( {x^p y^q } \right) + \Li_m (x^p )\Li_n (y^q )\,,
\end{split}
\end{equation}

where $\Li$ is a polylogarithm function.

\bigskip

Setting $p=0=q$ in equation~\eqref{equ.o6mnea9} or taking limit as $N\to\infty$ directly in equation~\eqref{equ.ybdhkk3} we have

\begin{equation}\label{equ.cfh5l7g}
\sum\limits_{r = 1}^\infty  {\frac{{H_{r,n} }}{{r^m }}}  + \sum\limits_{r = 1}^\infty  {\frac{{H_{r,m} }}{{r^n }}}  = \zeta (m + n) + \zeta (m)\zeta (n),\quad n,m\ne 1\,,
\end{equation}

The use of equations~\eqref{equ.m2dzb1g} and \eqref{equ.ybh1wcb} allows equation~\eqref{equ.cfh5l7g} to be written in terms of the polygamma function as

\[
\begin{split}
&\frac{{( - 1)^n }}{{(n - 1)!}}\sum\limits_{r = 1}^\infty  {\frac{{\psi _{n - 1} (r)}}{{r^m }}}  + \frac{{( - 1)^m }}{{(m - 1)!}}\sum\limits_{r = 1}^\infty  {\frac{{\psi _{m - 1} (r)}}{{r^n }}}\\ 
&\qquad\qquad = \zeta (m + n) + \zeta (m)\zeta (n),\quad n,m\ne 1\,.
\end{split}
\]

The particular case $m=n$ in equation~\eqref{equ.cfh5l7g} gives

\begin{equation}\label{equ.khvarp3}
2\sum\limits_{r = 1}^\infty  {\frac{{H_{r,n} }}{{r^n }}}  = \zeta (2n) + \zeta (n)^2,\quad n\ne 1\,.
\end{equation}

The result equation~(4.20) of reference~\cite{themis02} corresponds to an evaluation of the identity~\eqref{equ.khvarp3} at $n=2$.

\bigskip

Equation~\eqref{equ.khvarp3} is listed as Formula~(43) in~\cite{hnmathworld}.

\end{example}

\begin{example}

In the limit $N\to\infty$, identities~\eqref{equ.vse99ju} and \eqref{equ.r00mpzh} of Example~\ref{ex.e363h43} become

\begin{equation}\label{equ.oq0wqls}
\begin{split}
&\sum\limits_{r = 1}^\infty  {\frac{{h_{r,n} }}{{(2r - 1)^m }}}  + \sum\limits_{r = 1}^\infty  {\frac{{h_{r,m} }}{{(2r - 1)^n }}}\\
&\quad= (1 - 2^{ - m - n} )\zeta (m + n) + (1 - 2^{ - m} )(1 - 2^{ - n} )\zeta (m)\zeta (n),\quad n,m\ne 1\,,
\end{split}
\end{equation}

and

\begin{equation}
2\sum\limits_{r = 1}^\infty  {\frac{{h_{r,n} }}{{(2r - 1)^n }}}  = (1 - 2^{ - 2n} )\zeta (2n) + (1 - 2^{ - n} )^2 \zeta (n)^2,\quad n\ne 1\,,
\end{equation}

while identities~\eqref{equ.nrl5ek4} and \eqref{equ.hyn439v} become
\begin{equation}
\begin{split}
&\sum\limits_{r = 1}^\infty  {\frac{{h_{r,n} }}{{(2r + 1)^m }}}  + \sum\limits_{r = 1}^\infty  {\frac{{h_{r,m} }}{{(2r + 1)^n }}}  = \zeta (m)\zeta (n)(1 - 2^{ - m} )(1 - 2^{ - n} )\\
&\qquad\qquad\qquad\qquad\qquad\qquad\qquad\qquad - \zeta (m + n)(1 - 2^{ - m - n} )
\end{split}
\end{equation}

and

\begin{equation}\label{equ.ab0l0c2}
\sum\limits_{r = 1}^\infty  {\frac{{h_{r,n} }}{{(2r + 1)^m }}}  + \sum\limits_{r = 1}^\infty  {\frac{{h_{r,m} }}{{(2r - 1)^n }}}  = \zeta (m)\zeta (n)(1 - 2^{ - m} )(1 - 2^{ - n} )\,.
\end{equation}

In particular

\begin{equation}
2\sum\limits_{r = 1}^\infty  {\frac{{h_{r,n} }}{{(2r + 1)^n }}}  = \zeta (n)^2 (1 - 2^{ - n} )^2  - \zeta (2n)(1 - 2^{ - 2n} )
\end{equation}

and

\begin{equation}
\sum\limits_{r = 1}^\infty  {\frac{{h_{r,n} }}{{(2r + 1)^n }}}  + \sum\limits_{r = 1}^\infty  {\frac{{h_{r,n} }}{{(2r - 1)^n }}}  = \zeta (n)^2 (1 - 2^{ - n} )^2\,.
\end{equation}

\end{example}

\begin{example}
The result equation~(3.8c) of reference~\cite{zheng} implies that

\begin{equation}\label{equ.l1dnzjg}
\sum\limits_{r = 1}^\infty  {\frac{{h_{r,3} }}{{(2r - 1)^2 }}}  = \frac{{\pi ^2 }}{{16}}\zeta (3) + \frac{{31}}{{64}}\zeta (5)\,,
\end{equation}

from which, upon using identity~\eqref{equ.oq0wqls}, we get

\begin{equation}\label{equ.e9z9vtc}
\sum\limits_{r = 1}^\infty  {\frac{{h_{r,2} }}{{(2r - 1)^3 }}}  = \frac{{3\pi ^2 }}{{64}}\zeta (3) + \frac{{31}}{{64}}\zeta (5)\,.
\end{equation}

From identities~\eqref{equ.l1dnzjg} and \eqref{equ.e9z9vtc} and using identity~\eqref{equ.ab0l0c2} we get

\[
\sum\limits_{r = 1}^\infty  {\frac{{h_{r,3} }}{{(2r + 1)^2 }}}  = \frac{{\pi ^2 }}{{16}}\zeta (3) - \frac{{31}}{{64}}\zeta (5)
\]

and

\[
\sum\limits_{r = 1}^\infty  {\frac{{h_{r,2} }}{{(2r + 1)^3 }}}  = \frac{{3\pi ^2 }}{{64}}\zeta (3) - \frac{{31}}{{64}}\zeta (5)\,.
\]

Since

\[
\begin{split}
h_{_{r - 1,2} }^2  &= \left( {h_{r,2}  - \frac{1}{{(2r - 1)^2 }}} \right)^2\\
&= h_{r,2}^2  - \frac{{2h_{r,2} }}{{(2r - 1)^2 }} + \frac{1}{{(2r - 1)^4 }}\,,
\end{split}
\]
we also have

\[
2\sum\limits_{r = 1}^\infty  {\frac{{h_{r,2}^2 }}{{4r^2  - 1}}}  = 2\sum\limits_{r = 1}^\infty  {\frac{{h_{r,2} }}{{(2r - 1)^3 }}}  - (1 - 2^{ - 5} )\zeta (5)\,,
\]

from which we get, upon using equation~\eqref{equ.e9z9vtc}

\begin{equation}
\sum\limits_{r = 1}^\infty  {\frac{{h_{r,2}^2 }}{{4r^2  - 1}}}  = \frac{{3\pi ^2 }}{{64}}\zeta (3)\,.
\end{equation}

\end{example}

\begin{example}\label{ex.fomfjoe}

Letting $N\to\infty$ in identity~\eqref{equ.afc1dnp} of Example~\ref{ex.cu17le1} we obtain

\begin{equation}
\sum\limits_{r = 1}^\infty  {\frac{{h_{r,m} }}{{r^n }}}  + \sum\limits_{r = 1}^\infty  {\frac{{H_{r - 1,n} }}{{(2r - 1)^m }}}  = (1 - 2^{ - m} )\zeta (m)\zeta (n),\quad n\ne 1,\,m\ne 1\,.
\end{equation}

In particular,

\begin{equation}
\sum\limits_{r = 1}^\infty  {\frac{{h_{r,n} }}{{r^n }}}  + \sum\limits_{r = 1}^\infty  {\frac{{H_{r - 1,n} }}{{(2r - 1)^n }}}  = (1 - 2^{ - n} )\zeta (n)^2,\quad n\ne 1\,.
\end{equation}

\end{example}

\begin{example}

From the definition of $\bar H$ and the identities~\eqref{equ.to49yt7} and \eqref{equ.ldfid07} it follows that

\begin{equation}\label{equ.e2skunm}
\mathop {\lim }\limits_{N \to \infty } \bar H_{N,n}  = \left\{ \begin{array}{l}
 \log 2,\quad n=1 \\
 \\ 
 \left( {1 - \frac{1}{{2^{n - 1} }}} \right)\zeta (n),\quad n\ne 1 \\ 
 \end{array} \right.\,.
\end{equation}

Hence, from the identities~\eqref{equ.h87qnmg} and \eqref{equ.bx3kgp2} we obtain

\begin{equation}
2\sum\limits_{r = 1}^\infty  {( - 1)^{r - 1} H_{r} }  = \log 2
\end{equation}

and

\begin{equation}\label{equ.bhdgdgs}
2\sum\limits_{r = 1}^\infty  {( - 1)^{r - 1} H_{r,n} }  = \left( {1 - \frac{1}{{2^{n - 1} }}} \right)\zeta (n)\,,\quad n\ne 1\,.
\end{equation}

Similarly from identities~\eqref{equ.n7gklxg} and \eqref{equ.idbeoly} and using identities~\eqref{equ.gfrtfed} and \eqref{equ.nhgftfj} we have

\begin{equation}\label{equ.d6ajlzi}
\begin{split}
2\sum\limits_{r = 1}^\infty  {( - 1)^{r - 1} h_{r,n} }  &= \sum\limits_{r = 1}^\infty  {\frac{1}{{(4r - 3)^n }}}  - \sum\limits_{r = 1}^\infty  {\frac{1}{{(4r - 1)^n }}}\\
&= \frac{{( - 1)^{n} }}{{4^n \Gamma (n)}}\left\{ {\psi _{n - 1} \left( {\frac{1}{4}} \right) - \psi _{n - 1} \left( {\frac{3}{4}} \right)} \right\}\,.
\end{split}
\end{equation}

In reference~\cite{kolbig} it was established that
\[
\psi _{2n} \left( {\frac{1}{4}} \right) - \psi _{2n} \left( {\frac{3}{4}} \right) =  - \pi (2\pi )^{2n} |E_{2n} |
\]

and

\[
\psi _{2n - 1} \left( {\frac{1}{4}} \right) - \psi _{2n - 1} \left( {\frac{3}{4}} \right) = (2n - 1)!2^{4n} \beta (2n)\,,
\]

where

\[
\beta (m) = \mathop {\lim }\limits_{N \to \infty } \bar h_{N,m}  = \sum\limits_{s = 1}^\infty  {\frac{{( - 1)^{s - 1} }}{{(2s - 1)^m }}} 
\]
and $E_m$ is the $mth$ Euler number defined by the exponential generating function

\[
\frac{2}{{e^t  + e^{ - t} }} = \sum\limits_{m = 0}^\infty  {\frac{{E_m t^m }}{{m!}}}\,. 
\]

Using these results in identity~\eqref{equ.d6ajlzi} we obtain

\begin{equation}
2\sum\limits_{r = 1}^\infty  {( - 1)^{r - 1} h_{r,2n} }  = \beta (2n)
\end{equation}

and

\begin{equation}
2\sum\limits_{r = 1}^\infty  {( - 1)^{r - 1} h_{r,2n - 1} }  = \frac{{|E_{2n - 2}| }}{{2^{2n} \Gamma (2n - 1)}}\pi ^{2n - 1}\,.
\end{equation}

In particular

\begin{equation}
2\sum\limits_{r = 1}^\infty  {( - 1)^{r - 1} h_r }  = \frac{\pi }{4}\,,
\end{equation}

\begin{equation}
2\sum\limits_{r = 1}^\infty  {( - 1)^{r - 1} h_{r,2} }  = G
\end{equation}

and

\begin{equation}
2\sum\limits_{r = 1}^\infty  {( - 1)^{r - 1} h_{r,3} }  = \frac{{\pi ^3 }}{{32}}\,.
\end{equation}

From identities~\eqref{equ.n8tkvqw} and \eqref{equ.av4l176} we have

\begin{equation}\label{equ.cgzc8jh}
\sum\limits_{r = 1}^{\infty} {( - 1)^{r - 1} \left( {2H_{r,n} H_{r,m}  - \frac{{H_{r,n} }}{{r^m }} - \frac{{H_{r,m} }}{{r^n }}} \right)}=-\left ( 1-\frac{1}{2^{m+n-1}}\right )\zeta(m+n)\,.
\end{equation}

Setting $m=n$ in identity~\eqref{equ.cgzc8jh} yields

\begin{equation}
\begin{split}
2\sum\limits_{r = 1}^\infty  {( - 1)^{r - 1} \left( {H_{r,n}^2  - \frac{{H_{r,n} }}{{r^n }}} \right)}  &= 2\sum\limits_{r = 1}^\infty  {( - 1)^{r - 1} H_{r,n} H_{r - 1,n} }\\
&=  - \left( {1 - \frac{1}{{2^{2n - 1} }}} \right)\zeta (2n)\,.
\end{split}
\end{equation}

Thus

\begin{equation}\label{equ.xbjeju6}
2\sum\limits_{r = 1}^\infty  {( - 1)^{r - 1} H_{r,n} H_{r - 1,n} }  =  - \frac{{(2^{2n - 1}  - 1)}}{{(2n)!}}\left| {B_{2n} } \right|\pi ^{2n}\,,
\end{equation}

where $B_{m}$ is the $mth$ Bernoulli number defined by

\[
\frac{t}{{e^t  - 1}} = \sum\limits_{m = 0}^\infty  {B_m \frac{{t^m }}{{m!}}} \,.
\]

In particular

\begin{equation}\label{equ.dg360qb}
2\sum\limits_{r = 1}^\infty  {( - 1)^{r - 1} H_r H_{r - 1} }  = 2\sum\limits_{r = 1}^\infty  {( - 1)^{r - 1} H_r^2 }  - 2\sum\limits_{r = 1}^\infty  {( - 1)^{r - 1} \frac{{H_r }}{r}}  =  - \frac{{\pi ^2 }}{{12}}\,.
\end{equation}

From identity~\eqref{equ.dg360qb} and the known result
\[
2\sum\limits_{r = 1}^\infty  {( - 1)^{r - 1} \frac{{H_r }}{r}}  = \frac{{\pi ^2 }}{6} - \log ^2 2,\quad\mbox{(\cite{chu}, equation~4.2c)}
\]
we obtain

\begin{equation}
2\sum\limits_{r = 1}^\infty  {( - 1)^{r - 1} H_r^2 }  = \frac{{\pi ^2 }}{{12}} - \log ^2 2\,.
\end{equation}

Setting $m=0$ in identity~\eqref{equ.cgzc8jh} and using identities~\eqref{equ.e2skunm} and \eqref{equ.bhdgdgs} we obtain

\[
2\sum\limits_{r = 1}^\infty  {( - 1)^{r - 1} rH_{r,2} }  =  - \frac{{\pi ^2 }}{{24}} + \log 2
\]
and

\[
2\sum\limits_{r = 1}^\infty  {( - 1)^{r - 1} rH_{r,n} }  = \left( {1 - \frac{1}{{2^{n - 2} }}} \right)\zeta (n - 1) - \frac{1}{2}\left( {1 - \frac{1}{{2^{n - 1} }}} \right)\zeta (n)\,,\quad n\ne 1,n\ne 2 \,.
\]

From identities~\eqref{equ.gpqd08z} and \eqref{equ.v1b4oyo} we have

\[
2\sum\limits_{r = 1}^\infty  {( - 1)^{r - 1} h_{r,n} h_{r - 1,n} }  = \frac{1}{{4^{2n} \Gamma (2n)}}\left( {\psi _{2n - 1} \left( {\frac{3}{4}} \right) - \psi _{2n - 1} \left( {\frac{1}{4}} \right)} \right)=-\beta {(2n)}\,.
\]

In particular

\begin{equation}\label{equ.pq6lucu}
2\sum\limits_{r = 1}^\infty  {( - 1)^{r - 1} h_r h_{r - 1} }  =  - G\,.
\end{equation}

From the corrected version of equation~4.5c of~\cite{chu} 
\[
2\sum\limits_{r = 1}^\infty  {( - 1)^{r - 1} \frac{{h_r }}{{2r - 1}}}  = \frac{{\pi \log 2}}{4} + G
\]
and the identity~\eqref{equ.pq6lucu} we deduce that

\begin{equation}
2\sum\limits_{r = 1}^\infty  {( - 1)^{r - 1} h_r^2 }  = \frac{{\pi \log 2}}{4}\,.
\end{equation}

\end{example}

\section{Conclusion}

We have given and proved a summation identity which we subsequently applied in its various forms to obtain mostly new finite and infinite summation formulas involving the generalized harmonic numbers.

\medskip

\end{document}